\documentclass[11pt]{article}

\usepackage{graphicx}           
\usepackage{latexsym}
\usepackage{latexsym}
\usepackage{amsfonts}
\usepackage{amssymb}
\usepackage{amsmath}

\newcommand\vect[1]{{\bf#1}}
\newcommand\matr[1]{{\bf#1}}

\renewcommand{\tau}{c}

\newcommand{\beq}{\begin{equation}}
\newcommand{\eeq}{\end{equation}}
\newcommand{\be}{\begin{enumerate}}
\newcommand{\ee}{\end{enumerate}}
\newcommand{\bi}{\begin{itemize}}
\newcommand{\ei}{\end{itemize}}
\newcommand{\bc}{\begin{center}}
\newcommand{\ec}
{\end{center}}
\def\real{\hbox{\rm\setbox1=\hbox{I}\copy1\kern-.45\wd1 R}}

\begin{document}

\title
{Speeding-Up Convergence via Sequential Subspace Optimization: Current State and Future Directions
}
\author{Michael Zibulevsky \\  \\ Department of Computer Science\\ Technion---Israel Institute of Technology \\  E-mail: mzib@cs.technion.ac.il }
 \maketitle

\section*{Abstract}

This is an overview paper written in style of research proposal. In recent years we introduced a general framework for large-scale unconstrained optimization --  Sequential Subspace Optimization (SESOP) and  demonstrated its usefulness 
for  sparsity-based  signal/image denoising, deconvolution, compressive sensing, computed tomography,  diffraction imaging, support vector machines. 
We explored its combination with Parallel Coordinate Descent and Separable Surrogate Function methods, obtaining state of the art results in above-mentioned areas.

There are several  methods, that are faster than plain SESOP under specific conditions: Trust region Newton method - for  problems with easily invertible Hessian matrix;
Truncated Newton method -  when  fast multiplication by Hessian  is available;
Stochastic optimization methods - for problems with large stochastic-type data;
Multigrid methods - for problems with nested multilevel  structure. 
Each of these methods can be further improved by merge with SESOP.
One can also  accelerate Augmented Lagrangian method for constrained optimization problems  and Alternating Direction
Method of Multipliers for problems with separable objective function and non-separable constraints. 


 \section{Background}

 Many  problems in science and engineering are handled as optimization tasks, often with very high dimensions. 
Solving them calls for the use of iterative methods of various sorts, and then convergence speed becomes crucial.  
Interior-point  polynomial complexity methods  provide fast and robust solution of convex problems, where Newton optimization is applicable \cite{nester_nemir_book,  boyd_book}. 

When the problem size exceeds  several thousand of variables, storage and inversion of Hessian matrix required for the Newton method become prohibitively expensive. Therefore increased attention to the methods which use  gradient information only. Several of them possess optimal worst-case complexity \cite{Nemirovski-1983}: ORTH method by Nemirovsky  \cite{Nemirovski-1982}, Nesterov method \cite{Nesterov-1983, Nesterov-2003} and  FISTA method by Beck and Teboulle \cite{beck2009fast}.


Still  solving large ill-conditioned problems with high accuracy remains a challenge. Worst-case bounds are too pessimistic in many real-life cases.  Recently we have developed a technique  called SESOP, which constitutes a significant step towards  this goal.



\subsection{SESOP -- Sequential Subspace Optimization }
 The story of SESOP method \cite{Narkiss-2005} begins with the method of Conjugate Gradients (CG)
\cite{Stiefel-1952}. Quadratic CG (i.e. CG applied to a quadratic
function) has remarkable convergence properties: Its linear
convergence rate (see for example \cite{NocedWrite}) is
$\frac{\sqrt{r}-1}{\sqrt{r}+1}$,
where $r$ is the condition number of the Hessian of the objective. This rate is much
better than the steepest descent rate, $\frac{{r}-1}{{r}+1}$.

One can also rely on a $1/k^2$ sub-linear worst-case convergence of
the quadratic CG, which does not depend on the Hessian conditioning
(see for example \cite{Nemirovski-1982,Nemirovski-1983,Narkiss-2005}):
\begin{eqnarray}\label{eq:qadratic_decay}
f(\vect{x}_{k+1})- f_{\mbox{optimal}} \leq
\frac{L\Vert \vect{x}_0-\vect{x}_{\mbox{optimal}}\Vert ^2}{k^2}
\end{eqnarray}
where $k$ is the iteration index and  $L$ is the Lipschitz constant of the gradient
of $f$. The presented convergence rates are intimately related to the well known
{\em expanding manifold} property of quadratic CG: At every iteration the method
minimizes the objective function over an affine subspace spanned by directions of
all previous propagation steps and gradients.

In the case of a smooth convex function (not necessarily  quadratic), one could
propose a similar algorithm that preserves the expanding manifold property. Such an
algorithm should minimize the objective function over an affine subspace spanned by
directions of all previous propagation steps and the latest gradient. This method
inherits the $1/k^2$ convergence of CG,  however, the cost of an iteration of such a
method will increase with iteration count.

In order to alleviate this problem, Nemirovski \cite{Nemirovski-1982} suggested to
restrict the optimization subspace just to three directions: the current gradient,
the sum of all previous steps and a ``cleverly'' weighted sum of all previous
gradients (a version of such weights is given in \cite{Narkiss-2005}). The resulting
ORTH-method inherits the optimal worst case $1/k^2$ convergence
(\ref{eq:qadratic_decay}), but it does not coincide with CG, when the objective is
quadratic, and typically converges slower than CG, when the function become "almost"
quadratic in the neighborhood of solution.
%

The SESOP method \cite{Narkiss-2005} extends the ORTH subspaces with several directions of
the last propagation steps. This way, the method, while preserving  a $1/k^2$
convergence for smooth convex functions, becomes equivalent to the CG in the
quadratic case. This property boosts the efficiency of the algorithm.

 Quite often the function we minimize has a form $f(x)=\phi(Ax)$, where multiplication by matrix $A$ is costly, and computation of $\phi(\cdot)$ is relatively cheap.  The low-dimensional subspace optimization task at every iteration of SESOP can be
addressed using the Newton algorithm. The main computational burden in this process
is the need to multiply the spanning directions by $A$, but these multiplications
can be stored in previous iterations, thus enabling the SESOP speed-up with minor  additional cost.

 In order to improve efficiency further  we can substitute gradient direction with  direction of parallel coordinate descent (PCD)
or direction, provided by minimizer of a separable surrogate function (SSF)  \cite{MatalonJ} .
This approach provides state of the art results in the area of sparse signal approximation, where the objective function is 
\begin{eqnarray}\label{eq:BP}
f(\vect{x})=\| \matr{Ax} - \vect{b}\|_2^2 + \mu \| \vect{x}\|_1
\end{eqnarray}
\paragraph{ Parallel Coordinate Descent (PCD)} Quite often coordinate descent is faster than gradient descent method in terms of total computational burden, because re-evaluating function value may be very cheap while changing one coordinate. In particular, for function \eqref{eq:BP} coordinate optimization can be calculated analytically and  involves just one column of matrix $\matr{A}$.
Still for large dense matrices of size $n\times n$ we need $2n^2$ operations for one pass over all coordinates, i.e. the cost of two matrix-vector multiplications.

 On the other hand in many problems fast matrix-vector multiplication is available, e.g. Fast Fourier Transform or Fast Wavelet Transform. For such problems we use Parallel Coordinate Descent \cite{ShrinkElad}:
Staying at current point we evaluate coordinate descent steps for all coordinates without moving along them (this can be performed analytically at  cost of two fast matrix-vector multiplications). Then we obtain the next iterate of $\vect{x}$ via approximate minimization of the objective function (line search) along  the  obtained vector of coordinate steps.

\paragraph{Method of separable surrogate functions (SSF)} is another efficient alternative to gradient descent  proposed by \cite{Daub, Figu1}. At every iteration 
the first term of  \eqref{eq:BP} is substituted with diagonal quadratic function which has the same value  as this term  at current iterate and majorates it elsewhere. Combining this diagonal quadratic term with $\mu \| \vect{x}\|_1$, we obtain a majorating separable surrogate function (SSF) for our problem. At every iteration we  minimize SSF analytically, which always provides reduction of $f(\vect{x})$. Similarly to PCD,  SSF step  costs  about two matrix-vector multiplications, and the speed of convergence of two methods  is also comparable.  Still when the problem is ill-conditioned, the convergence may be slow.  

\paragraph{PCD-SESOP and SSF-SESOP method}  We can further accelerate the convergence if at every iteration we minimize the objective function over affine subspace spanned by current PCD or SSF direction and several previous propagation steps. 
PCD-SESOP and SSF-SESOP methods  
are asymptotically equivalent to  the diagonally preconditioned CG. On the other hand, the PCD and
SSF directions provide much faster progress at initial steps, when compared to the
ordinary nonlinear CG. This partially explains extremal efficiency of these methods
on difficult problems, where they are often significantly faster  \cite{ZibEladSPM} than popular Nesterov  method and FISTA \cite{nesterov2007gradient, beck2009fast}.

\begin{figure}
    \begin{center}
        \begin{tabular}{cc}
            \includegraphics[width=0.47\columnwidth, height=0.45\columnwidth]{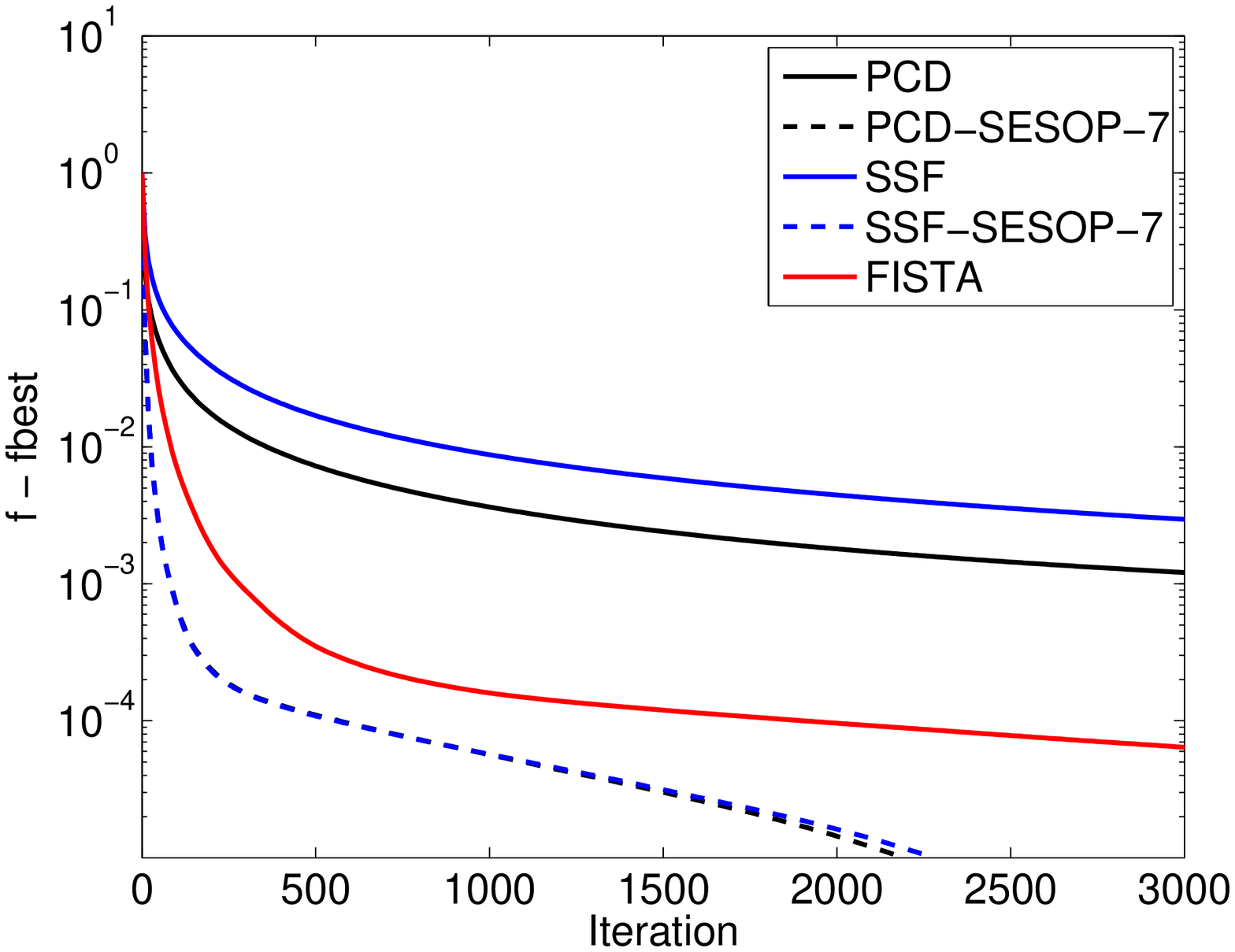}&
            \includegraphics[width=0.47\columnwidth, height=0.45\columnwidth]{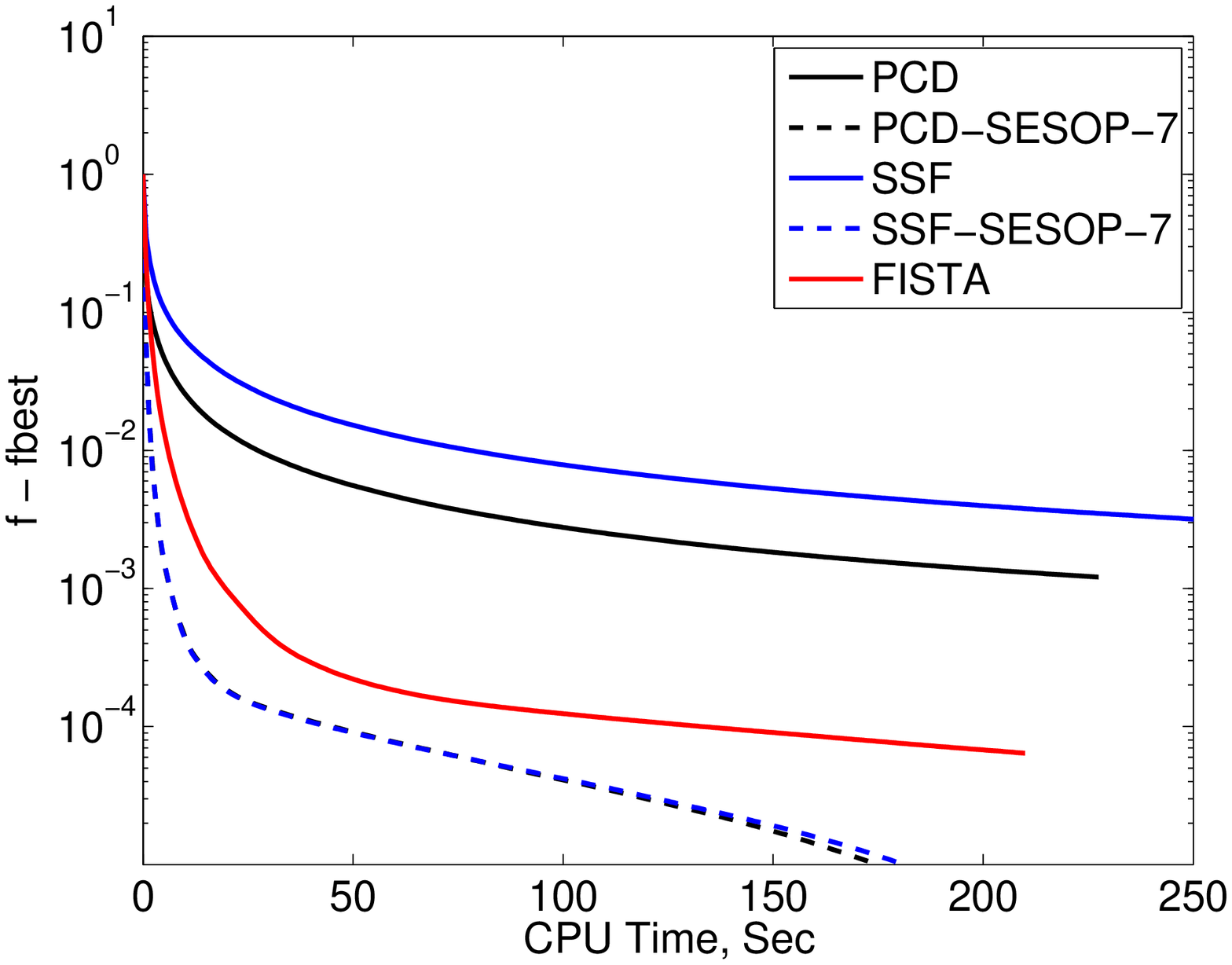}\\
            \includegraphics[width=0.47\columnwidth, height=0.45\columnwidth]{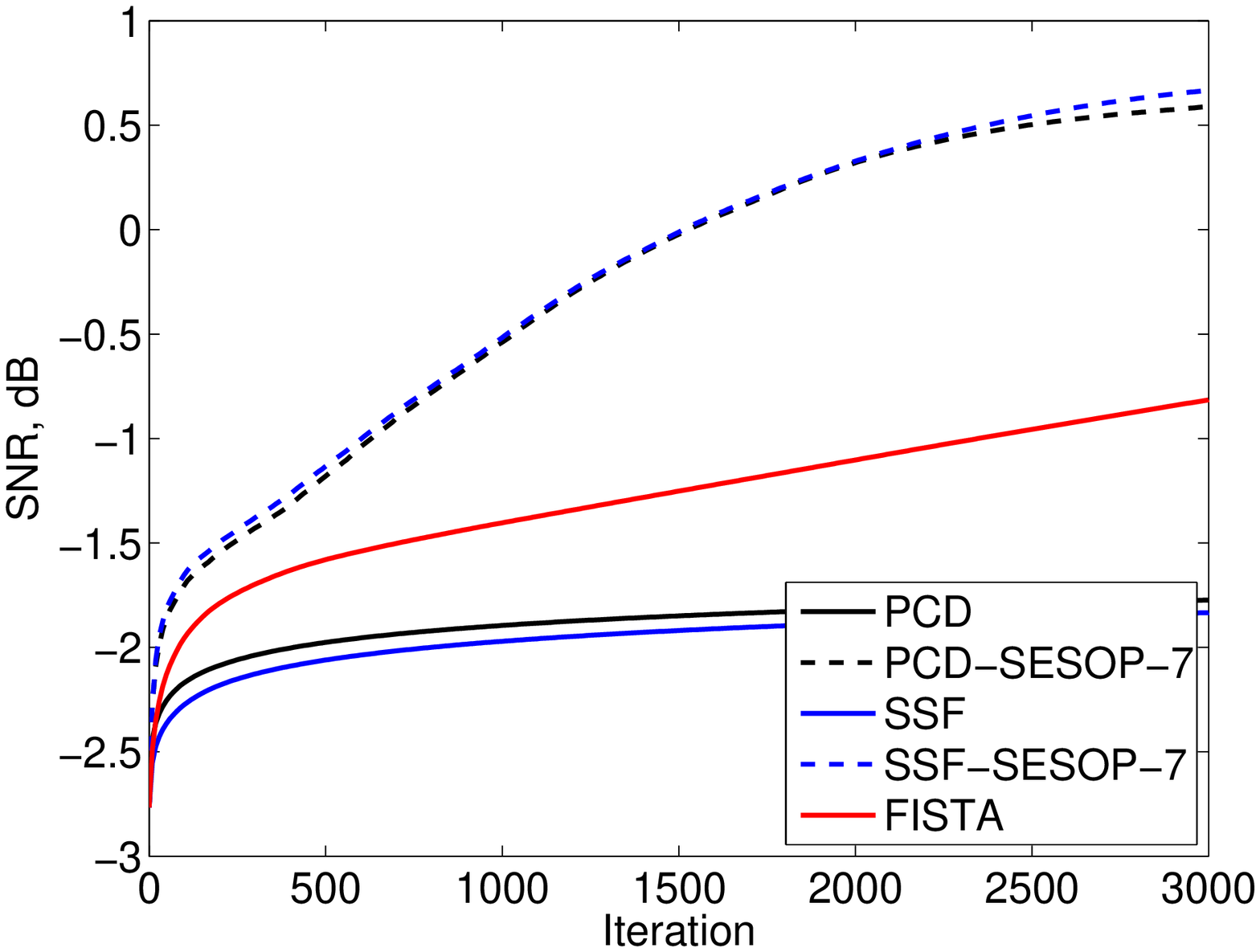}&
            \includegraphics[width=0.47\columnwidth, height=0.45\columnwidth]{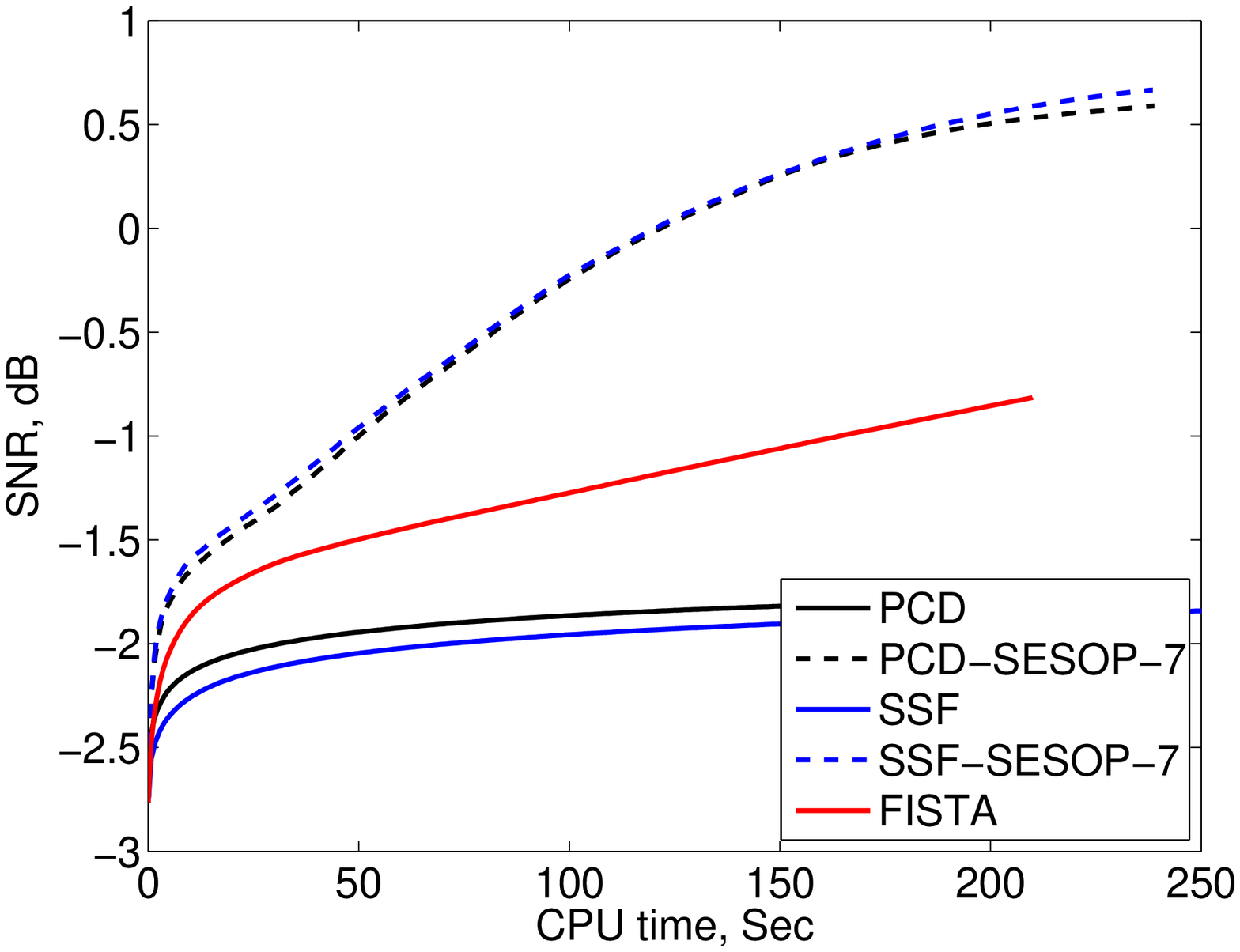}
        \end{tabular}
    \end{center}
    \caption{\label{fig:Ignace2} Comparison of optimization methods over ill-posed geo-science type problem  with matrix $\matr{A}$
of size $1848 \times 8192$}
\end{figure}

In Figure \ref{fig:Ignace2} we present results of an experiment  \cite{ZibEladSPM} adopted
from \cite{Ignace}.  The problem uses an explicit matrix $\matr{A}$
of size $1848 \times 8192$, termed $K_4$, which is built from
a Gaussian random matrix, by shaping its singular values to fit a
classical and highly ill-conditioned geo-science problem  
Beyond the natural goal of comparing various algorithms, in this
experiment we aim to address two additional issues: (i) Run-time and
its relation to number of iterations; and (ii) The behavior of the
algorithms was explored  for very small value of $\mu=1e-6$, for which
PCD/SSF methods are known to deteriorate. Figure
\ref{fig:Ignace2} presents the objective function and the SNR of recovery of $x$, both
as a function of the iteration count, and as a function of time.
As we can see,  FISTA converges much faster then SSF, which
corresponds to the observations in \cite{Ignace}. On the other hand,
PCD-SESOP and SSF-SESOP are far superior to all other methods.

\subsection{Newton-type  optimization methods}

As a part of our proposal we are going to merge SESOP with several well-known Newton methods in order to boost their performance. In this section we just set up  the  notions.
There are two basic approaches in continuous multidimensional unconstrained optimization: trust region and line search. Each of them has its own strength and weakness. 
Suppose that we are looking for a minimum of a function of several variables
\beq\label{eq1}
\min ⁡f(x)
\eeq
If we have an easily   minimizable  approximate model  $q_0 (x)$  of our function around the initial point $x_0$, we could find an approximate minimizer of \eqref{eq1} as
$x_1= \arg  \min⁡〖q_0 (x)〗$
In Newton method, for example, $q_0(\cdot)$ is a second-order Taylor expansion of $f(\cdot)$ around $x_0$. We can continue with this process iteratively in order to improve accuracy of solution
$
x_{k+1}= \arg  \min⁡〖q_k (x)〗
$
The main problem with this approach is that our model $q_k (x)$ may become inaccurate far from its base point $x_k$ , therefore function value in the next point $x_{k+1}$may even increase. In order to alleviate this problem, one can use line search or trust region strategy.

\paragraph{Line search method}
At iteration k we first compute a minimizer of the model
\beq
y_{k}= \arg \min⁡〖q_k (x)〗
\label{eq4}
\eeq
and then perform one-dimensional optimization in its direction 
\beq
d_k=y_{k}-x_k 
\label{eq5}
\eeq
\paragraph{Trust region method}
A potential disadvantage of line search based optimization is that for ill-behaved functions,  even though the model $q$ fits the function quite well locally, the global minimum of the model may be far away from the minimum of the original function, and direction $d$ may be "poor". Therefore an alternative trust region approach restricts search of the optimum of the model by a limited area around the current iterate $x_k$
\begin{eqnarray} 
x_{k+1}&=\arg \min q_k(x)  \label{eq_trust_region}\\ 
\mbox{subject to:} &||x-x_k||\leq r_k  \nonumber 
\end{eqnarray} 
Parameter $r_k$ is adjusted dynamically from iteration to iteration based on progress in actual function value, see e.g. \cite{NocedWrite}.
The problem with trust region approach is two fold. First, finding constrained minimum of the model is often 2 -- 3 times more expensive than finding the unconstrained minimum. Secondly, the model can fit the original function well in some directions and poorly in others, therefore an ideal trust region should be a kind of ellipsoid instead of Euclidean ball. However adjusting parameters of such ellipsoid is rather difficult. In our proposal we will alleviate this difficulty using SESOP.

\paragraph{Truncated Newton (TN) method} \cite{Dembo-1983, Nash-2000} is used when computing and inverting Hessian matrix is prohibitively expensive. Therefore at every outer iteration we minimize  quadratic Taylor expansion  $q_k(x)$  around the current iterate~$x^k$ only approximately, using limited number of CG steps. 
The outer iteration of TN is accomplished with a line search,
in order to guarantee function decrease.
The overall effectiveness of the TN method is rather sensitive to
the choice of stopping rule for the internal CG optimization. We attempt to overcome this
difficulty, replacing  the line search with subspace optimization.
In this way we allow the CG iterations to stay matched through consequent
TN steps.

\subsection{Mini-batch stochastic optimization}
Machine learning poses data-driven optimization problems in which the objective function involves the summation of loss terms over a set of data to be modeled. Classical optimization techniques must compute this sum in its entirety for each evaluation of the objective, respectively its gradient. As available data sets grow ever larger, such "batch" optimizers therefore become increasingly inefficient. They are also ill-suited for the online (incremental) setting, where partial data must be modeled as it arrives. Stochastic (online) gradient-based methods, by contrast, work with gradient estimates obtained from small subsamples (mini-batches) of training data. This can greatly reduce computational requirements: on large, redundant data sets, simple stochastic gradient descent routinely outperforms sophisticated second-order batch methods by orders of magnitude in spite of the slow convergence of first-order gradient descent.

Recently there were partially successful attempts to adapt classical powerful unconstrained optimization methods to mini-batch mode. For example of mini-batch limited memory quasi-Newton L-BFGS  \cite{StochQN, CarefulQN, OptMethDeepLearning2011}. It significantly outperforms stochastic gradient method, however still there is  room for improvement.

\section{Merging  existing algorithms with with SESOP}

The methods mentioned above are faster than plain SESOP under specific conditions: Trust region Newton method -- for  problems with easily invertible Hessian matrix;
Truncated Newton method --  when  fast multiplication by Hessian  is available;
Stochastic optimization methods - for problems with large stochastic-type data;
Multigrid methods -- for problems with nested multilevel  structure.  Each of these methods has its own weakness, which  may be  alleviated using subspace optimization.

Another task is to accelerate Augmented Lagrangian method for constrained optimization problems  and Alternating Direction
Method of Multipliers for problems with separable objective function and non-separable constraints, using SESOP concept.


\subsection{Merging trust region and line search approaches within SESOP framework}
The known problem with trust region approach \eqref{eq_trust_region} is two fold. First, finding constrained minimum of the model is often 2 -- 3 times more expensive than finding the unconstrained minimum. Secondly, the model can fit the original function well in some directions and poorly in others, therefore an ideal trust region should be a kind of ellipsoid instead of Euclidean ball. However adjusting parameters of such ellipsoid is rather difficult. In our proposal we will alleviate this difficulty using SESOP.

Usually the model $q_k$ in 
trust region step \eqref{eq_trust_region} is built in the way that its gradient at the point $x_k$ is equal to the gradient of the original function $f$. This is the case, for example, for Newton method, when the model is the second order Taylor expansion of $f$ around $x_k$. Therefore, if the radius $r_k$ of trust region goes to zero in \eqref{eq_trust_region}, the obtained direction will correspond to the one of steepest descent. On the other hand, when the radius goes to infinity, the direction will correspond to the pure Newton step. For intermediate values of the radius the obtained direction gives a compromise between Newton and gradient step.

Motivated by the above observation, we propose to compute the next iterate  $x_{k+1}$ as a minimizer of the original function $f$ over affine subspace that touches $x_k$ and spans the current "Newton" direction $d_k$ of the line search method given by \eqref{eq4}, \eqref{eq5} and the direction of negative gradient $\ -\nabla f(x_k)$. In addition it may be very useful to include  directions of several previous steps $x_k - x_{k-1}$, $\ x_{k-1}-x_{k-2} $, ..., as well as previous gradients and "Newton" directions into the subspace. 

The proposed method takes advantage of several worlds. Whenever the Newton direction is nor effective, it is at least as efficient as SESOP \cite{Narkiss-2005,  Narkiss-2005, MatalonJ, ZibEladSPM}, which is usually much faster than steepest descent or nonlinear conjugate gradient method; when Newton direction is good, local quadratic convergence rate of Newton  is preserved.

We should also note that quite often subspace optimization may be performed in very efficient way, so that additional computational burden is very small comparing to the evaluation of Newton direction. It is possible to incorporate the proposed approach in other newton-type schemes, for example into Levenberg-Marquardt method for nonlinear least squares, which is very popular e.g. in training of moderately-sized feed-forward neural networks.

\subsection{Combining  SESOP with Truncated Newton method}

Like we mentioned before, TN method is  sensitive to  change in early stopping of  the internal CG optimization. We are going to resolve this problem 
 replacing  the line search in the outer step with subspace optimization. In this way the CG sequence will not break between outer iterations and  will  stay matched through consequent
TN steps.

When minimizing objective function $f(x)$, $k$-th TN step will approximately minimize its  quadratic model around current iterate $x^k$
\beq
q_k(x)=f(x^k)+{g^k}^T(x-x^k) + \frac{1}{2}(x-x^k)^T H_k (x-x^k), \label{QuadModel}
\eeq
where $g^k=\nabla f(x^k)$ is the gradient and $H_k=\nabla^2 f(x^k)$ -- the
Hessian of $f$ at $x^k$. 
Suppose that  TN step $k$ was truncated after $l$ CG iterations. Coming
back from the quadratic model 
to the original objective, 
we would like to imitate the next CG step, using
$f(x)$ instead of $q(x)$. CG iteration $l+1$ inside TN would perform optimization of
the quadratic model $q(x)$ over the affine subspace $S_{k,l}$, passing through the current inner
iterate $x^{k,l}$ and  spanned by the last CG step $x^{k,l}- x^{k,l-1}$ and the current
gradient $\nabla q(x^{k,l})$. Instead, the next SESOP iteration will minimize $f$
over the extended affine subspace $S_k \supset S_{k,l}$.
In order to provide monotone descent of $f$, we add to $S_k$  the TN direction
$$d_{TN}=x^{k,l}-x^k.$$ Now $x^k\in S_k$, and any monotone method used for the subspace
optimization over $S_k$ starting from $x^k$,  will reduce the objective relatively
to $f(x_k)$.
Optionally, we include several previous outer steps and gradients of $f$ into $S_k$,
in order to improve the function descent, when the TN directions are not good enough.

\paragraph{Next Truncated Newton step} After performing the subspace optimization,
we start a new TN iteration. At this stage, in order to keep alignment through the global CG
sequence, we perform the first new CG step as an optimization of the new quadratic model
$q_{k+1}(x)$  over the 2D subspace spanned by  $ x^{k+1}-
x^{k,l}$ and $g(x_{k+1})$.

\paragraph{\underline{Summary of SESOP-TN algorithm:  outer iteration $k$}}
\begin{enumerate}

\item  {\bf TN step}\ \ \ Solve approximately Newton system  $\nabla^2 f(x^k)d^k_{TN}=-\nabla f(x^k)$,
i.e. minimize quadratic model $q_k(x)$ in \eqref{QuadModel}, using $l$ steps of
CG. Denote the last CG iterate as $x^{k,l}$.

\item {\bf Subspace optimization step} $ \ \ \ x^{k+1} \approx \arg\min_{x\in S_k} f(x)$,

 where
affine subspace $S_k$ passes through $x^{k}$ and is spanned by:

   $\bullet$ TN Direction $d^k_{TN} = x^{k,l}-x^k$;\\
     $\bullet$  Last value of the gradient of quadratic model $\nabla q_k(x^{k,l})$  used in TN; \\
     $\bullet$ Last used CG direction in TN: $(x^{k,l}-x^{k,l-1})$;\\
     $\bullet$ {\small \nolinebreak{[Optional] directions of several previous outer steps and gradients of $f$.}}
\item  {\bf Goto TN step},
while performing the first new CG step as an optimization of quadratic model
$q_{k+1}(x)$ over 2D subspace spanned by  $ (x^{k+1}- x^{k,l})$ and $\nabla
f(x^{k+1})$.

\end{enumerate}
The presented procedure resolves the problem of TN sensitivity to  early break of the CG
process. For example, when the objective function is quadratic, SESOP-TN trajectory coincides
with the trajectory  of CG  applied directly to the objective function,
independently of the stopping rule in the TN step (see Figure \ref{quad_fig} in Section \ref{TN_simulations} below)
Standard TN lacks this property and converges more slowly when truncated too early.

\subsection{Mini-batch stochastic optimization using SESOP}

Inspired by the recent success  of mini-batch  limited memory quasi-Newton method L-BFGS  \cite{StochQN, CarefulQN, OptMethDeepLearning2011}  in various machine learning tasks,  we plan to adapt SESOP  to similar conditions. SESOP has shown state of the art efficiency on many classes of large scale problems, see e.g  \cite{MatalonJ, ZibEladSPM}. Basic SESOP iteration consists in computing minimum of the objective function over subspace spanned by the current gradient and several previous steps. In stochastic mode one can perform optimization of partial objective function, computed in  current mini-batch, over  subspace spanned by directions of steps in several previous  mini-batches and the previous partial gradient. Restricting use of the latest gradient we hope to stabilize the method by preventing over-fitting to the latest mini-batch data.

\subsection{Incorporating SESOP  into multigrid/multilevel optimization}

In  multigrid and multilevel techniques (see basic introduction in \cite{IradWhyMultigrid})  the problem whose solution is sought is represented in a
hierarchy of resolutions, and each such version of the problem is ``responsible'' for eliminating errors on a scale compatible with its typical resolution. 
 Multigrid methods accelerate computation in two ways. First, an approximate solution to the fine-grid problem may be cheaply obtained using coarse-grid formulation. Second, by "clever" sequential iteration between fine and coarse grid, an accurate solution can be obtained much faster, then just iterating on fine grid with a good initial guess \cite{IradWhyMultigrid}.

\paragraph{Merging SESOP with  multigrid} 
Our approach is relevant to both, linear and non-linear multigrid problems, expressed in variational form  as nested set of quadratic or general non-linear optimization problems, see e.g. \cite{nash2000multigrid,Toint_TrustRegion_Multiscale}.
Quite often optimization at every level of resolution  is carried by Conjugate Gradient (CG) method, which converges  much faster than steepest-descent type methods on ill-conditioned problems. One the other hand, CG efficiency is greatly reduced when the method is restarted after small number of iterations (restarts happen because of interleave of fine and coarse grid steps in multigrid.) In order to alleviate this problem, we propose to use SESOP  instead of CG. It is known that SESOP is equivalent to CG when used for minimization of quadratic function \cite{Narkiss-2005}, but usually converges faster than CG, when the function is non-quadratic.

Standard iteration of SESOP consists of optimization of the objective function over affine subspace spanned by directions of several previous steps, current [preconditioned] gradient (or other reasonable directions, like parallel coordinate descent,  \cite{MatalonJ}). Additionally we are going to include  descent direction, provided by most recent coarse grid solution, into current fine-grid SESOP subspace. This should provide best combination of SESOP and multigrid advantages,  because fine-grid iterations have no break, and CG-type acceleration properties are preserved.

 \paragraph{Saddle point problem and Alternating Direction
Method of Multipliers} As a general framework, SESOP allows acceleration of many other memoryless techniques. In particular, we are going to explore saddle point algorithms such as  Augmented Lagrangian method for constrained optimization. In this case at every iteration of SESOP we will look for a saddle point in subspace of previous steps in primal and dual variables and current Augmented Lagrangian step.  When the primal objective is separable and Alternating Direction
Method of Multipliers is applicable (see e.g.  \cite{Boyd_alternating_2011}), we can include corresponding alternating directions  into SESOP subspace.

\section{Perspectives of  convergence analysis}
In this section we just share several thoughts regarding possible convergence analysis of SESOP related methods.
\paragraph{Composite functions} As we already mentioned, similarly to ORTH and Nesterov method, SESOP possesses optimal $1/k^2$ worst-case convergence expressed in  \eqref{eq:qadratic_decay}. 
The error in this formula is proportional to the Lipschitz constant $L$ of the gradient of the objective function $f$.  This bound become unsatisfactory when dealing with composite objective of type \eqref{eq:BP}, which incorporates non-smooth L1-norm term. Nesterov method for composite objective and FISTA alleviate this difficulty using SSF direction instead of the gradient in the update formula for evaluating  the next step. Therefore only Lipschitz constant of the gradient of the first smooth term in the objective is involved in error estimate.

In the same spirit SSF-SESOP method uses SSF direction and several previous steps when computing the next step. In all our numerical tests SSF-SESOP outperformed FISTA in number of iterations.  All the above rises hope that it is possible to develop a proof of optimal worst-case complexity of SSF-SESOP, which is similar to FISTA. 

On the other hand, asymptotic linear convergence rate of SESOP is similar to the Conjugate Gradient method and is much better than the Nesterov one. 
We have some results for quadratic functions \cite{MatalonJ}, however they could be extended for smooth nonlinear case as well as to the composite non-smooth functions.

\paragraph{Constrained optimization with PCD-SSF-SESOP}  Substituting $L_1$-norm  in \eqref{eq:BP} with $L_\infty$, and applying PCD-SSF-SESOP, one can obtain efficient method for problem with box-constraints. In similar way we can use indicator function of other simple feasible sets. The convergence results in this case should just follow from general results for composite functions.

\section{Some preliminary experiments: \\Truncated Newton SESOP (SESOP-TN)}\label{TN_simulations}
\begin{figure}
\begin{center}
\begin{tabular}{cc}
SESOP-TN & Standard TN\\
\includegraphics[width=0.47 \textwidth]{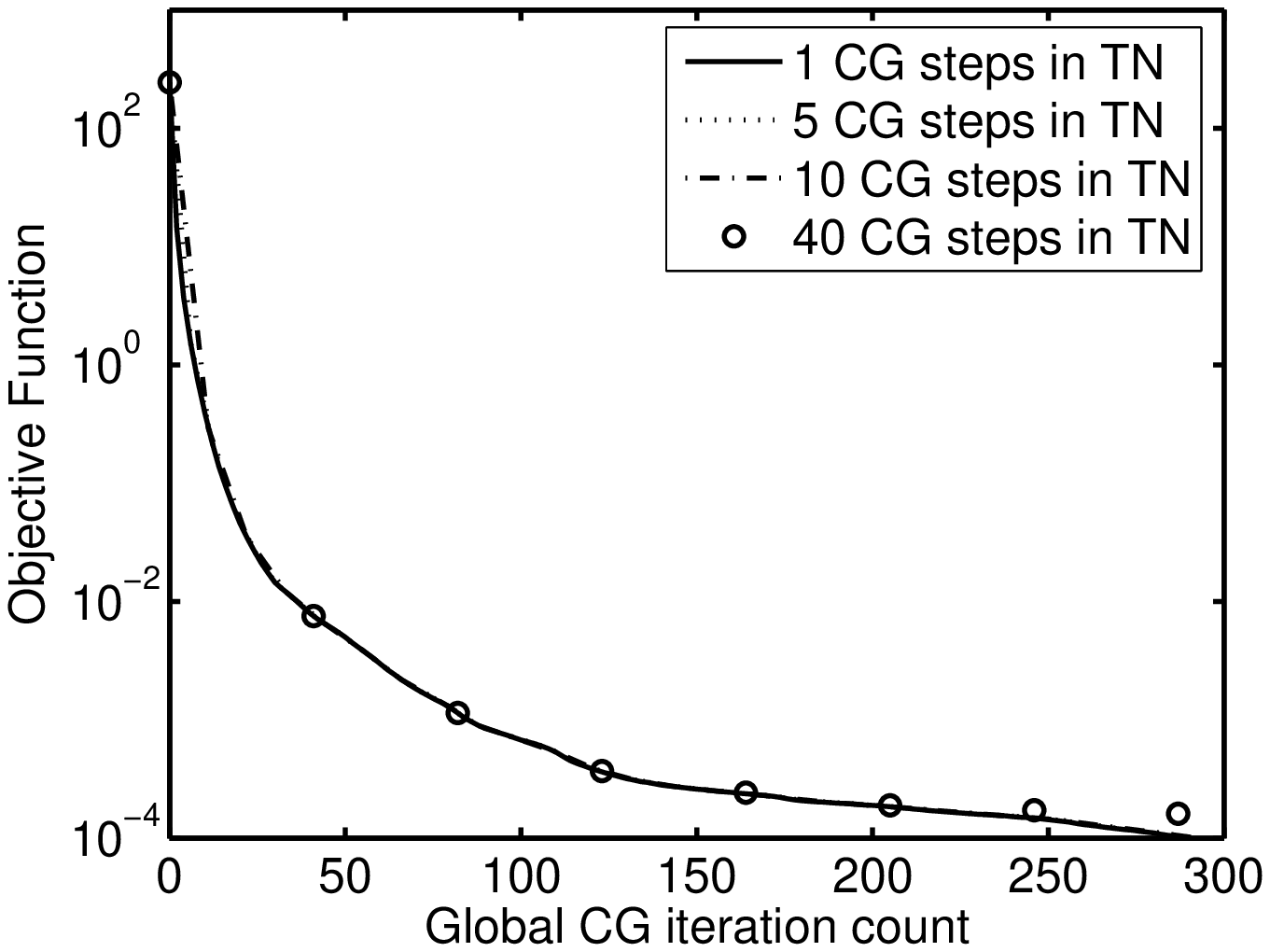} &
\includegraphics[width=0.47 \textwidth]{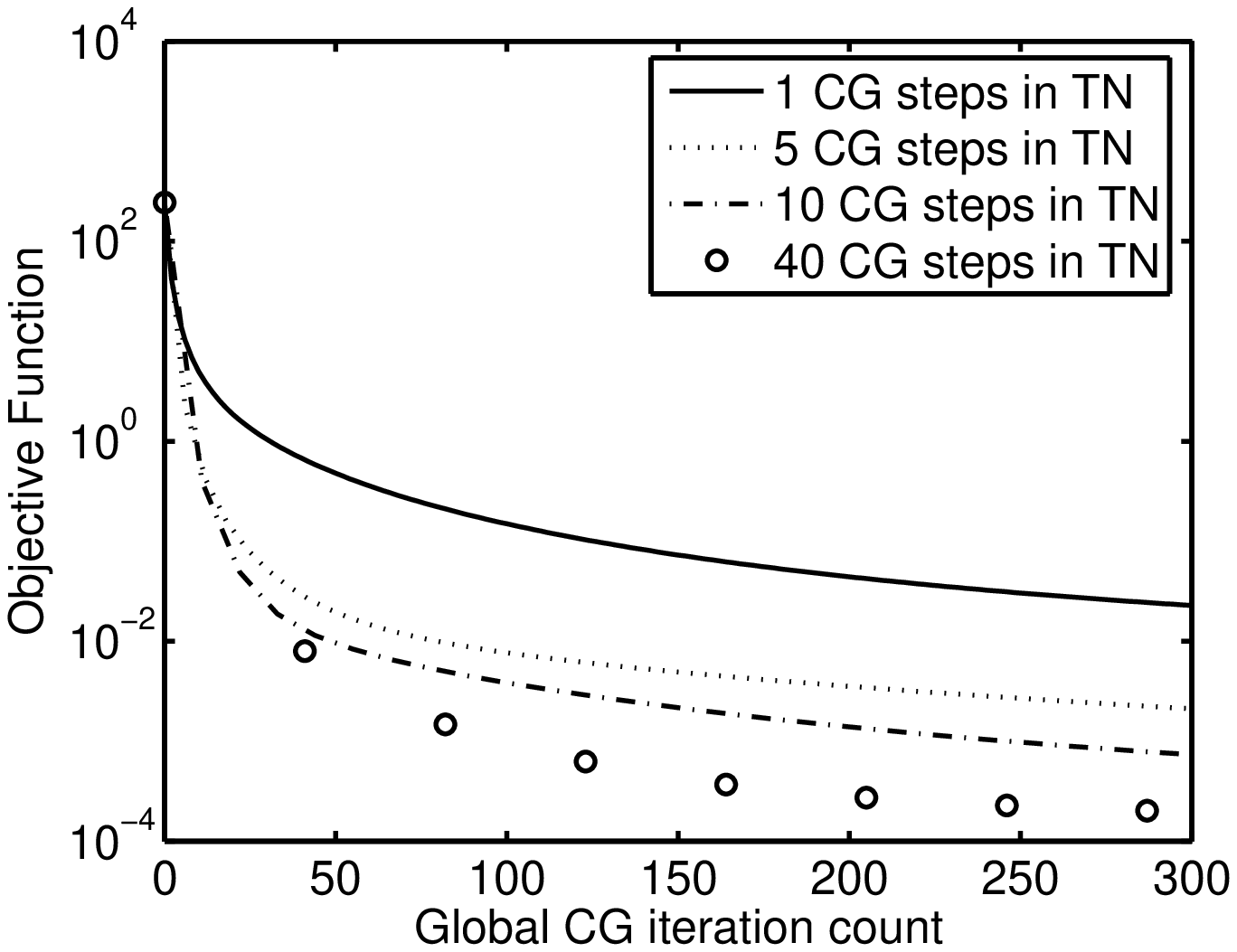}
\end{tabular}
\caption{Solving linear least squares \eqref{Eq-LS}, with  400 variables. The SESOP-TN
trajectory  does not depend on  the number of CG iterations in TN step.
Standard TN converges more slowly, when CG is truncated too early.
} \label{quad_fig}
\end{center}
\end{figure}

\paragraph{Quadratic function} First, let us demonstrate "proof of the concept" using a pure quadratic function.
We solve the linear least squares problem

\begin{equation}
\min||Ax-b||^2 
\label{Eq-LS}
\end{equation}
with $n=400$ variables, where the square random matrix  $A$ has
zero-mean {\em i.i.d.} Gaussian entries with variance  $1/n$.
As we see in Figure \ref{quad_fig}, SESOP-TN trajectory, as expected,  does not depend on the number of CG
iterations in the TN step. Standard TN (the right plot) lacks this property.

\begin{figure}
\begin{center}
{\em Exponents-and-Squares}, $200$ variables\\
\includegraphics[width=0.8 \textwidth]{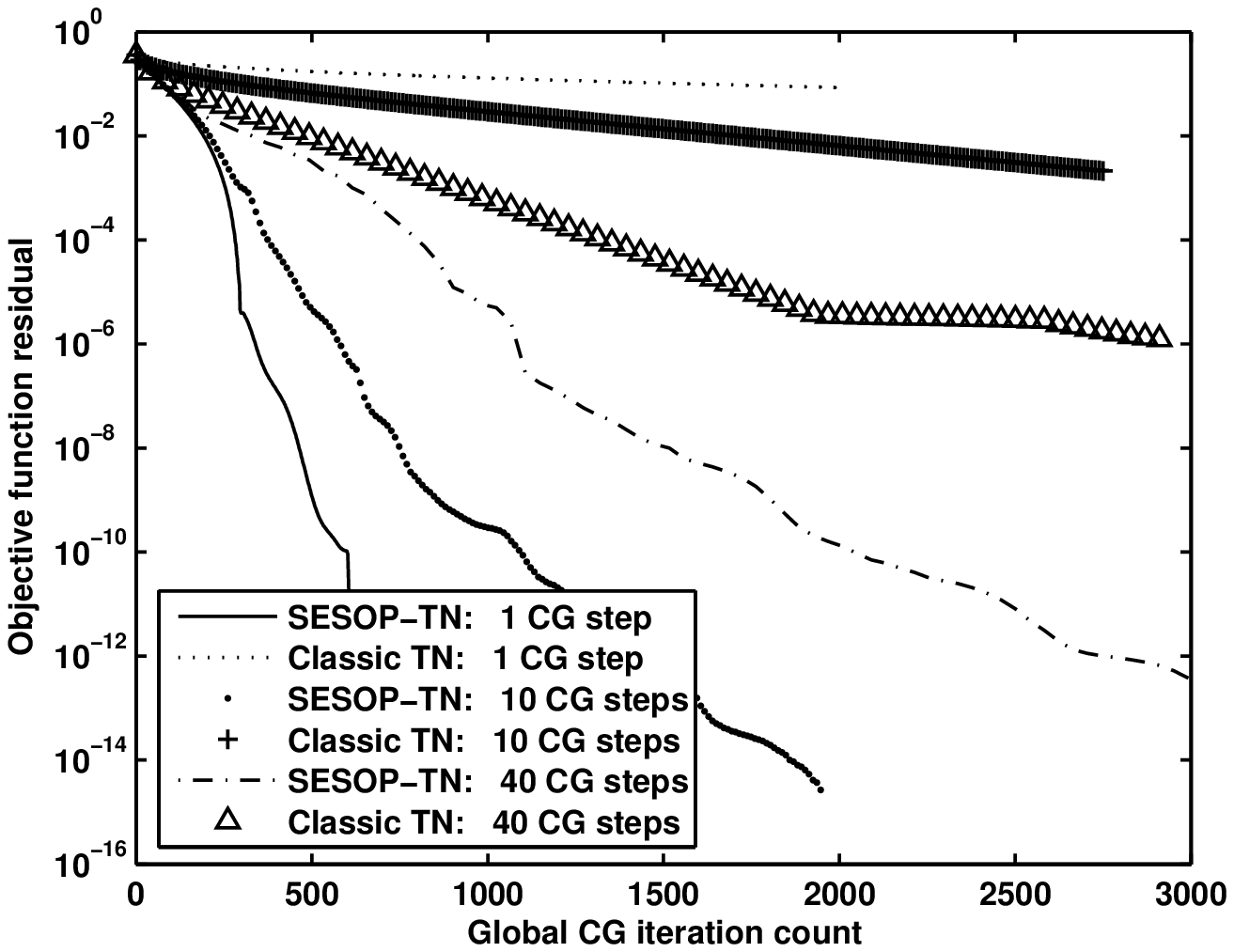}\\ \ \\

Linear SVM, 99758 variables\\
\includegraphics[width=0.8 \textwidth, height=0.6 \textwidth]{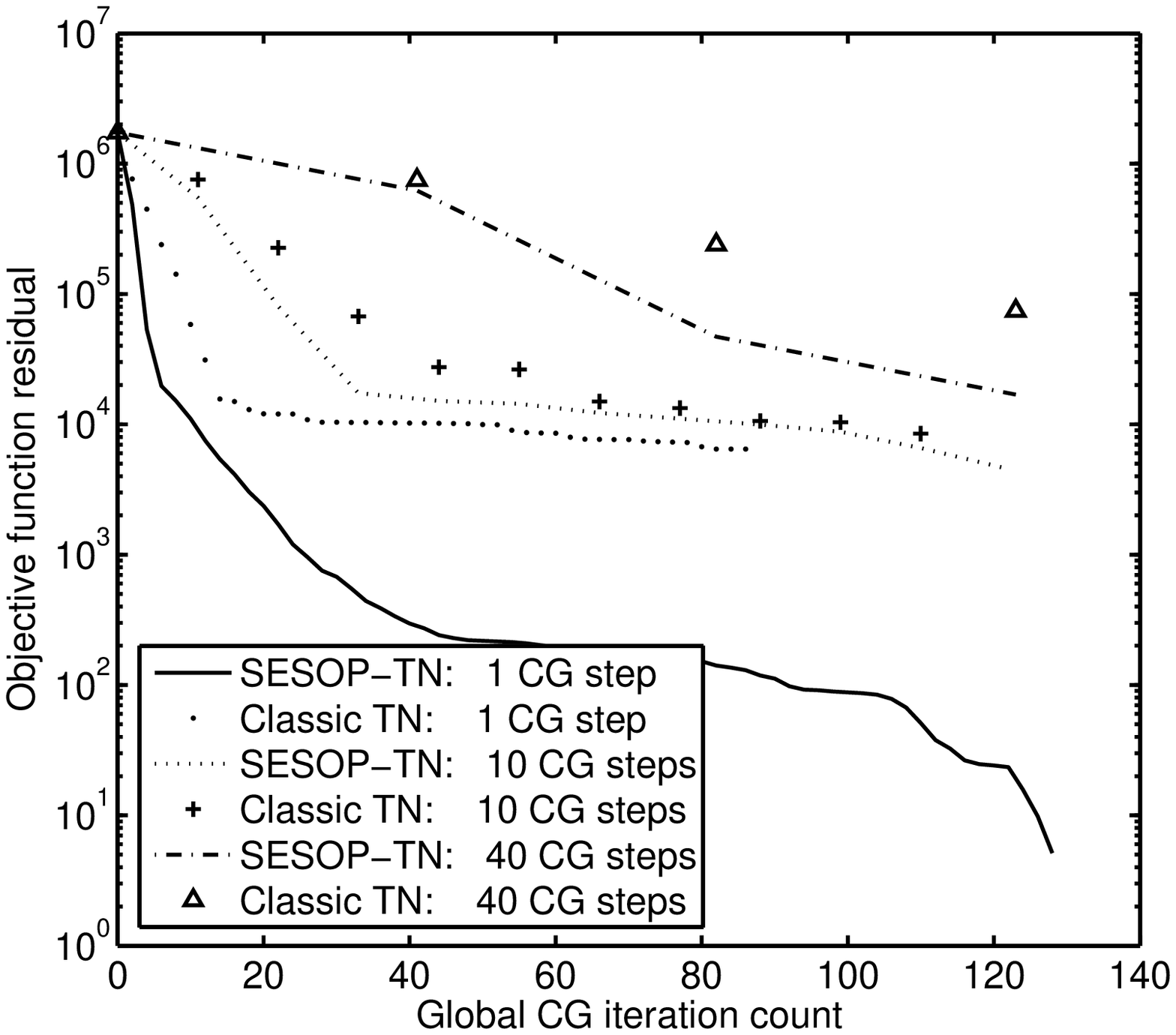}
\end{center}
\caption{Two nonlinear problems. The plots show the residual between the current objective and the optimal value
versus CG iteration count.
} \label{nonlin_fig}
\end{figure}

\paragraph{Two non-linear examples}   The first  problem is
{\em Exponents-and-Squares} \cite{Neilsen-2000} with $n=200$ variables:
$$ \min \ \ e^{-{\bf 1}^T x}+
\frac{1}{2}\sum_{j=1}^n j^2 x_j^2\ .$$
The second example is {\em Linear Support Vector Machine (SVM)},  see
\cite{Narkiss-SVM-2005} for more details on unconstrained formulation of SVM.
We used  data set {\em Astro-physics-29882} \cite{Joachims-2006} with 99758 variables,
and selected randomly 1495  training examples from there.
In both problems (see Figure~\ref{nonlin_fig}), SESOP-TN consistently outperformed classic TN, when  restricted
to 1, 10 or 40 CG iterations in TN step.

%
%
%

\section{Conclusions}
In this paper we presented several perspective directions of acceleration of optimization algorithms using SESOP framework and shared some thoughts about  convergence analysis.


\bibliography{sesop_proposal_ISFa}

\end{document}